\documentclass{elsarticle}
\usepackage{amssymb}
\usepackage[english]{babel}
\usepackage[letterpaper,top=2cm,bottom=2cm,left=3cm,right=3cm,marginparwidth=1.75cm]{geometry}
\usepackage{listings}
\usepackage{amsmath}
\usepackage{graphicx}
\usepackage{url}
\usepackage{xcolor}
\usepackage[colorlinks=true, allcolors=blue]{hyperref}
\usepackage{comment}

\newcommand{\qmod}[0]{$Q_{mod}$}
        
\begin{document}
\begin{frontmatter}
\title{MapperEEG: A Topological Approach to Brain State Clustering in EEG Recordings}

\author[inst1]{Brittany Story \corref{cor}}
\author[inst2,inst4]{Zhibin Zhou}
\author[inst2]{Ramesh Srinivasan}
\author[inst1]{Scott Kerick}
\author[inst1]{David Boothe}
\author[inst1,inst3]{Piotr J.~Franaszczuk}

\affiliation[inst1]{organization={U.S.~Army DEVCOM Army Research Laboratory},
addressline={7101 Mulberry Point Rd}, 
city={Aberdeen Proving Ground},
postcode={21005}, 
state={MD},
country={USA}}

\affiliation[inst2]{organization={University of California, Irvine},
addressline={2201 Social \& Behavioral Sciences Gateway}, 
city={Irvine},
postcode={92697}, 
state={CA},
country={USA}}

\affiliation[inst3]{organization={Johns Hopkins University School of Medicine},
addressline={733 N Broadway}, 
city={Baltimore},
postcode={21205}, 
state={MD},
country={USA}}

\affiliation[inst4]{organization={University of California, Davis Medical Center},
addressline={4150 V Street, Room 1200}, 
city={Sacramento},
postcode={95817}, 
state={CA},
country={USA}}

\cortext[cor]{Corresponding author}

\begin{abstract}
\noindent\textbf{Background}\\
Topological data analysis (TDA) has exploded as a tool for analyzing and making sense of high dimensional datasets across a variety of fields. Mapper is a tool from TDA that captures low-dimensional structure from high-dimensional data, precisely the approach needed to capture relevant information from high-dimensional neural time series. Electrical potential scalp recording, or electroencephalography (EEG), is routinely used in clinical applications and research studies thanks to its noninvasive nature, relatively inexpensive equipment, and high temporal resolution. But, it is prone to contamination, exhibits low spatial resolution, and has a non-stationary nature. Thus, it requires advanced signal processing and mathematical analysis methods for tasks requiring unsupervised brain state clustering. 

\noindent\textbf{New Method}\\
We introduce MapperEEG, an approach to unsupervised brain state clustering that uses tools from classical EEG analysis combined with Mapper to cluster and connect brain states.

\noindent\textbf{Results}\\
We show that MapperEEG can serve as a clustering algorithm in the spectral domain and provide additional information about the underlying brain state connectivity in a tapping task. Additionally, we use a go/no-go shooting task to explore how MapperEEG can still provide insight into the underlying structure and clusters of brain states even when it and other clustering methods fail.

\noindent\textbf{Comparison with Existing Methods}\\
We demonstrate that it outperforms six other clustering algorithms such as hierarchical clustering, Hidden Markov Models, and basic autoencoders on identifying states in a tapping task.

\noindent\textbf{Conclusions}\\
M-EEG offers a novel and effective approach to analyzing EEG data, showing promise for brain state clustering and analysis.
\end{abstract}

\end{frontmatter}

\section{Introduction}\label{sec:Intro}

Characterizing normal and pathological brain states via electroencephalography (EEG) has long been a key method in neuroscience.
EEG data is easy to collect, inexpensive, non-invasive, and provides high time resolution, which makes it a useful tool for researchers across a wide range of disciplines. 
It has some drawbacks, such as low spatial resolution and high susceptibility to noise and artifacts. 
Various analytical methods have been used for EEG data interpretation, both in clinical and research applications \cite{Kaplan2007HandbookOE}.
From the earliest EEG studies, clinicians and researchers noticed the presence of oscillatory features in EEG recordings. 
Over time, the oscillatory components of EEG within different frequency bands were correlated with particular brain functions and pathologies \cite{buzsaki2006rhythms}.
Other brain data collection methods, such as fMRI, were introduced later with higher spatial resolution and less noise, however, these methods provide lower time resolution and significantly higher costs. 

Recently, topological approaches have been used by computational neuroscientists to capture and cluster brain states from fMRI data. 
In \cite{saggar2018towards}, Saggar et al.~studied brain state organization using an algorithm called Mapper, introduced by Singh et al.~in \cite{singh2007mapper}. 
Mapper, when applied to fMRI data, was able to organize different brain states categorically without any pre-labeling. 
But, the types of tasks and states that researchers can study with this method are limited to those that can be collected in laboratory settings, subject to the inherent restrictions imposed by fMRI. 
Thus, we look to see how we can combine Mapper with tools from traditional EEG approaches to cluster EEG data to study the underlying structure of brain states. 

In this work, we introduce MapperEEG, an approach for analyzing EEG data that uses the power spectral density paired with the Mapper algorithm to construct a clustering algorithm for EEG states.
We demonstrate how combining these approaches can be used to understand brain state structure by showing its efficacy at clustering brain states during a human teaming tapping task and an executive control task under low and high stress conditions.
The paper is organized as follows. 
First, we review a variety of works that use different approaches for clustering EEG data in Section~\ref{sec:RelatedWork}.
Section~\ref{sec:Background} outlines relevant background information necessary for establishing MapperEEG.
Then, we outline the two datasets used for comparison in Section~\ref{sec:Datasets}.
Section~\ref{sec:Methods} describes the MapperEEG algorithm and explains the method step-by-step. 
Next, Section~\ref{sec:Results} compares MapperEEG's performance to other clustering approaches on the two datasets and explains some of the benefits provided by MapperEEG as opposed to other methods.
Then, in Section~\ref{sec:Discussion} we discuss our findings, the importance of the spectral domain, and current limitations. 
Finally, we outline conclusions, plans for future work, and extensions in Section~\ref{sec:Conclusion}.

\section{Related Work}\label{sec:RelatedWork}

Several approaches have investigated clustering brain states via EEG signals.
Common approaches, such as k-means, Hidden Markov Models, and microstate analysis have shown promising results on certain tasks like identifying concealed information~\cite{bablani2020multi} identifying eyes open/eyes closed~\cite{Jajcay2023micro}, and depth of anesthesia (DOA) monitoring~\cite{si2021mirco}.
As in the works cited above, combining EEG preprocessing techniques with clustering methods often yields more robust clustering algorithms that are more resilient to noise and artifacts. 
For example, Redwan et al.~\cite{Redwan2024psd} used signal power derived from the power spectral density (PSD) combined with several different machine learning algorithms, like random-forest models, support vector machines, and Gaussian process classifiers, to cluster brain states. 
Thus, we draw inspiration from approaches like Redwan et al.'s to combine traditional EEG preprocessing approaches with novel clustering algorithms to form a new model.

Several other tools have been used to cluster brain states.
Graphs have long been used to capture functional connectivity between brain states, identify state changes due to external stimuli, and assist with neurological disease identification~\cite{yu2018acupuncture, yu2019acupuncture, yu2020network}.
Different machine learning approaches have been used to classify brain states from EEG data, such as CNNs~\cite{Lawhern2018-ie} and support vector machines~\cite{Moctezuma2019svm, ieracitano2019svm}.
Saeidi et al.~\cite{Saeidi2021-jr} provides a comprehensive review of using machine learning for neural decoding of EEG signals.
Recent advances use graphs combined with machine learning to classify different EEG states \cite{Belhadi2025-gu, Almohammadi2023revealing}.
More recently, topological methods have been used for analyzing EEG data \cite{Altindis2021-ph, Xu2021-cu}.
These works support our approach of using topological and graph-based methods combined with traditional preprocessing to construct MapperEEG.

\section{Background}\label{sec:Background}
Below, we outline necessary background on EEG processing, topological data analysis, and graph metrics. 

\subsection{EEG}\label{subsec:EEG}

Electroencephalography (EEG) has been a popular neuroimaging technique for studying brain activity since its first human application in 1924 \cite{StLouis2016EEG}. 
As mentioned in the introduction, some advantages of EEG include its non-invasive and inexpensive data collection, along with its high-temporal resolution~\cite{he2018EEG}. 
Modern systems can use anywhere from 4 to 256 electrodes, and years of research have pinpointed the optimal electrode placement for a wide range of tasks, such as neural speech tracking, epilepsy, and concentration \cite{montoya2021speech, maher2023seizure, siamaknejad2019concentration}.
But using EEG for brain state clustering has disadvantages, just as with any data collection approach.
Lack of spatial resolution, individual brain and skull differences, depth restrictions, and noise are just a few of the difficulties researchers face \cite{hochberg2006sensors}.
Over the years, researchers have worked around these obstacles by using tools that clean and parse EEG data. 
Source localization, band separation, filtering, downsampling, power spectral density (PSD), and other signal processing techniques have allowed researchers to pull meaningful patterns and information from the raw EEG data \cite{hu2019eeg}.  

Understanding how brain states differ across tasks has been explored across a variety of domains such as during transitions, resting, seizures, and cognitive tasks \cite{kringelbach2020brainstates, Lee2012clus, Centeno2014-mj, finn2017funcon}.
Often, these approaches rely on other types of neuroimaging, such as fMRI, which can be limited to laboratory-restricted tasks \cite{glover2011fmri}.
Thus, there remains a need for a methodology that can cluster different brain states from EEG data across a wide range of tasks. 
Traditional methods such as event-related potentials and oscillatory dynamics have long been used to analyze and interpret EEG data \cite{buzsaki2006rhythms}.
More modern approaches such as Bayesian networks, support vector machines (SVMs), and random forests, characterize signals and associate them to behavior. 
For example, Bird et al.~have used these methods to classify brain states as neutral, relaxed, or concentrating with an accuracy of roughly 87\% \cite{bird2018testing}. 
EEG states have also been used as a method of communication for people with communication issues \cite{lazarou2018interface}.
Additionally, EEG has been used to classify teammate trust with a decision tree classifier with an accuracy level of 72\% \cite{firoz2022trust}.
Our goal is to develop an EEG clustering/mapping algorithm that achieves high accuracy while effectively capturing the relationships between brain states. 

\subsection{Simplicial Complexes and Persistent Homology}\label{subsec:PH}

Topology is a branch of mathematics that focuses on the shape of spaces. 
One way to study shape is to look at features that do not change under continuous deformations, such as scaling, twisting, or bending. 
Under these deformations, the number of components, holes, and voids do not change \cite{Hatcher2000alg_top}. 
Thus, counting the number of these features is one way to analyze the shape of a given space.
But, data points have no inherent interesting shape, they float in the underlying state space.
Thus, if we want to study the shape of a dataset, we have to transform the data into an object where we can calculate topological features, such as the number of components or holes a dataset contains.
One approach is to turn our dataset into a simplicial complex.

A \textbf{simplicial complex}~\cite{Hatcher2000alg_top} is a set of points, line segments, triangles, tetrahedra, and their n-dimensional counterparts, called \textbf{simplices}, with the following additional structure: 
\begin{enumerate}
    \item for any face in the simplicial complex, its sub-faces must also be in the simplicial complex, and
    \item any non-empty intersection of faces is also in the simplicial complex.
\end{enumerate}

It is possible to calculate the shape of simplicial complex through \textbf{simplicial homology}, a tool from topological data analysis that provides counts of the number of components, holes, voids, and their higher dimensional counterparts.
Thus, simplicial homology can be used to characterize structure and to understand the dataset's ``shape''.
There are multiple ways to convert datasets into simplicial complexes, but for this work, we will focus on Mapper.

\subsection{Mapper}\label{subsec:Mapper}

The Mapper algorithm was first developed by Singh, Memoli, and Carlsson in 2007 as a way to visualize and understand the structure of high dimensional data \cite{singh2007mapper}. 
Since its inception, it has been used in several different domains, such as fraud detection, gerrymandering, and neuroscience \cite{mitra2021fraud, thatcher2021cartographic, hasegan2024neuro}.
With respect to neuroscience, Mapper has been used to detect and identify transition states, pinpoint causal states, and predict future states \cite{geniesse2019dyneusr, tassi2024gastro, vannoni2024state}.
Other topological tools, such as persistent homology, have been used to look at simulated EEG data \cite{nasrin2019top}, cancer research \cite{loughrey2021cancer}, and brain injuries \cite{Nielson2015brain_injury}. 
These applications all leverage topology's ability to analyze data in a way that is structurally relevant.
In this work, we will focus on using Mapper to examine the underlying structure of data collected from multiple EEG channels and how these states are related across different tasks and over time.

\begin{figure}[h]
    \centering
    \includegraphics[width=0.99\linewidth]{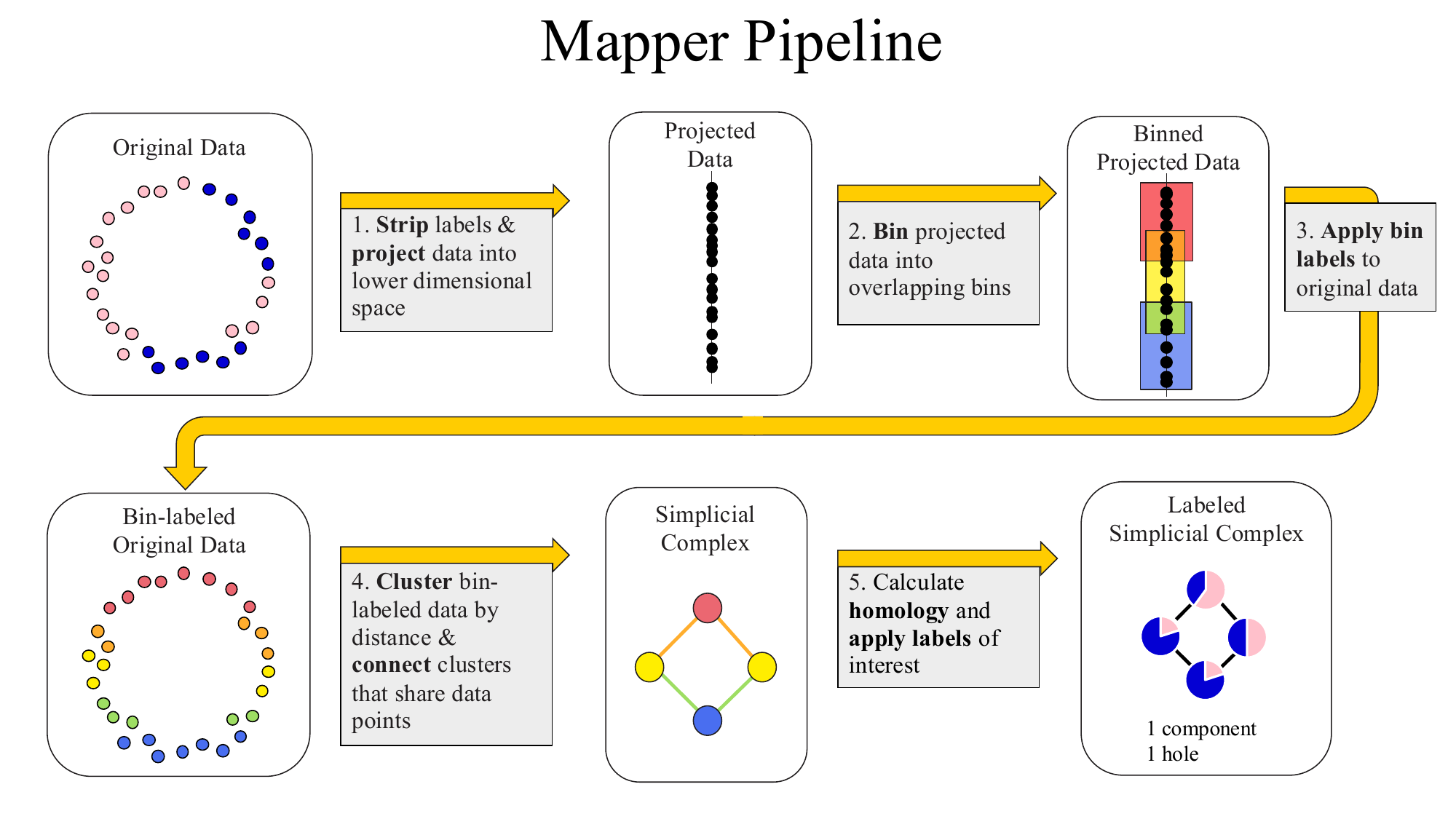}
    \caption{An example of the Mapper algorithm applied to a toy dataset. Note, the structure of the Mapper graph captures the shape of the dataset.}
    \label{fig:ex_pipeline}
\end{figure}

Mapper allows us to capture the structure of a dataset by turning the dataset into a simplicial complex and calculating its homology.
To build understanding, we outline a visualizable example for how the Mapper algorithm works. 
Consider a set of data points in $\mathbb{R}^2$ that roughly fall in the shape of a circle, see Figure~\ref{fig:ex_pipeline}.
\begin{enumerate}
    \item Project the label-stripped points into a lower dimension. For this example, let $f:\mathbb{R}^2\rightarrow\mathbb{R}^1$ be the projection function such that $f(x,y)=y$.
    \item Bin the projected points into $b$ bins with a $ov$\% overlap. For this example, let $b=3$ and $ov=33\%$. 
    \item Apply the bin labels of the projected points to the original points. Note, a point can have multiple bin labels, which we denote with orange and green in this example. 
    \item Cluster the points into nodes by their bin labels and a distance metric. Connect two nodes if they share at least one data point, building a simplicial complex. 
    \item Apply the original labels to each node, characterizing the distribution of the data points that were clustered into each node. In the example below, each node is a pie chart denoting the make up of labels of the data points grouped in that node. 
\end{enumerate}

Figure~\ref{fig:ex_pipeline}, outlines the Mapper pipeline, which captures the structure of the circular shape of the dataset with a simplicial complex with one component and one hole. 

For this paper, we will consider simplicial complexes consisting only of vertices and edges and we will ignore higher dimensional features such as triangles and solid tetrahedra. 
With these restrictions in place, we can refer to our simplicial complexes as \textbf{graphs}.
Further, we will restrict our focus to the number of components in a graph, ignoring higher dimensional topological features like holes or voids.
With that language in mind we will refer to the simplicial complex output of the Mapper algorithm as a \textbf{Mapper graph}.

\subsection{Q-mod Function}\label{subsec:qmod}

We seek to identify which Mapper graphs best capture the community structure of the data.
Particularly, we use modularity, specifically the \qmod function, as a way to evaluate Mapper graphs.
Modularity measures how well data are partitioned into communities where data points with the same label are connected within the same community when compared to graphs of the same size connected randomly. 
Garcia et al.~summarize how modularity has provided valuable insight into task performance, memory performance, and working memory \cite{garcia2018community}.
We use the $Q_{mod}$ function mentioned in Newman et al.~\cite{newman2004community} and Fortunato et al.~\cite{fortunato2010community}, similar to how it was applied in \cite{saggar2018towards}.

Let $G=\{V,E\}$ be the Mapper graph outlined above where $V$ is the set of vertices and $E$ is the set of edges. 
Thus, we define \textbf{\qmod} as, 
\[Q_{mod} = \frac{1}{2m}\sum\limits_{ij \in E}\left[A_{ij}-\frac{k_i k_j}{2m}\right]\delta(C_i, C_j),\]
where $A_{ij}$ is the adjacency matrix entry for row $i$ and column $j$, and $m$ is the number of edges in the graph.
For our purposes, $k_i$ denotes the strength of node $v_i$, which we choose to be the degree of $v_i$.
The community of the node, $C_i$, is the majority of labels within a node.
Note, labels can be pulled from the trial type, such as synchronized vs syncopated or low stress vs high stress.
If the labels are evenly split, we randomly select a label.
Finally, the function $\delta(C_i, C_j)$ is 1 when the communities of the vertices match and 0 otherwise.

The \qmod function captures the modularity of the graph by assigning a value between $(-\infty, 1)$.
A value approaching 1 represents a modular graph where vertices with the same labels are connected and there are few to no connections between vertices with different labels.
A graph with \qmod $\approx$ 0 captures the community structure of a graph where labeled vertices are connected randomly, some with the same labels and some with different labels.
A graph with \qmod $<0$ has fewer edges between vertices with the same labels and/or vertices connected to vertices with different labels \cite{fortunato2010community}.

\section{Datasets}\label{sec:Datasets}

We evaluate our approach on two datasets; a tapping task where participants are instructed to tap together over a variety of conditions and a shooting task where subjects are instructed to shoot targets in high and low stress conditions.
For visualizations of the paradigms, see Figure~\ref{fig:paradigms} and Figure~\ref{fig:paradigms_go}.
Below, we outline both datasets and define how we use the power spectral density (PSD) to construct the power data used within our MapperEEG algorithm.

\subsection{Tapping Task}

Neuroscientists have used synchronized/syncopated tapping tasks in a variety of contexts. 
Mayville et al.~used synchronized and syncopated tapping to study the cortical and subcortical brain networks via fMRI data \cite{Mayville2002ss}.
Additionally, Mayville et al.~have also studied how event-related changes in brain activity captured by magnetoencephalography (MEG) are associated with synchronous and syncopated tapping \cite{mayville2001ss_beta}.
But, both of these studies rely on fMRI or MEG data and as such, we want to show that EEG can also be used to differentiate between synchronization and syncopation. 

\begin{figure}[ht]
    \centering
    \includegraphics[width=0.999\linewidth]{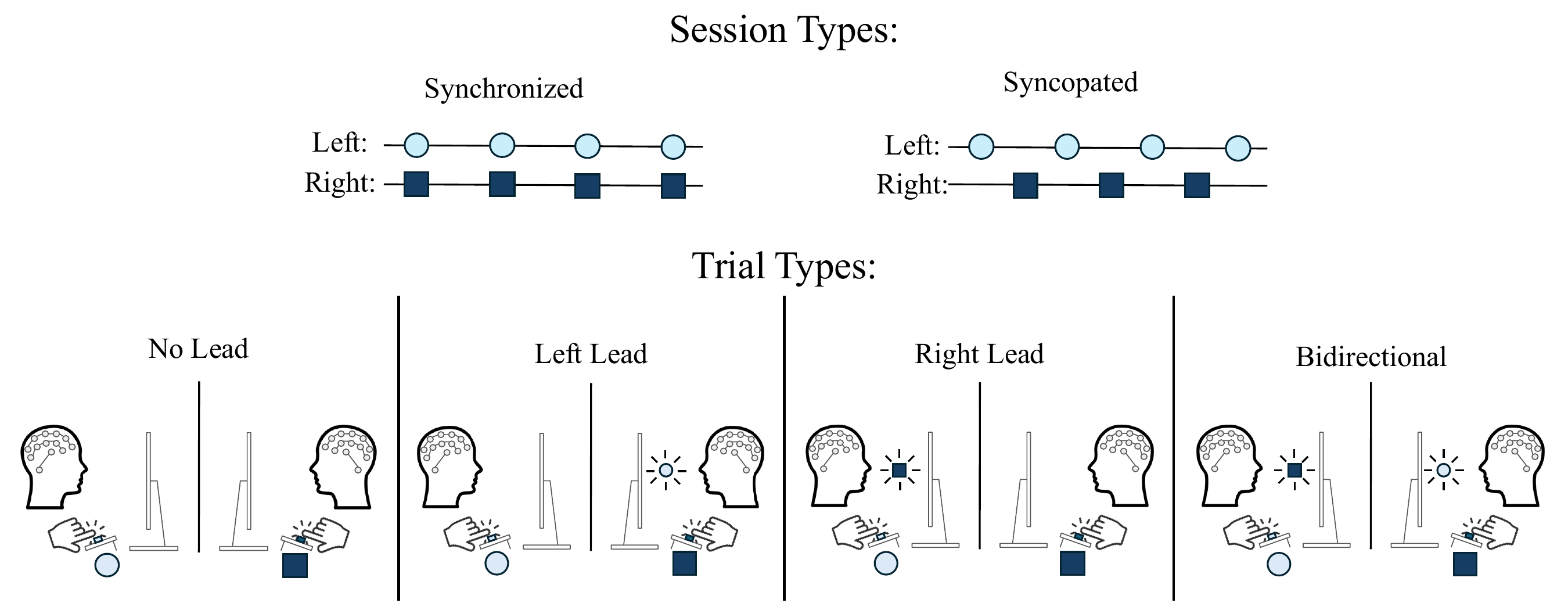}     
    \caption{The tappers were isolated from each other and participated in two sessions: a synchronized session where the tappers matched the beat and a syncopated session where the tappers tapped between the beat.
    During each session, the pairs experienced three runs of each of the four trial types: no lead, left lead, right lead, and bidirectional.
    No lead: both tappers tapped independently, with no feedback from the other tapper. 
    Left lead: the left tapper was instructed to synchronize/syncopate with the metronome and then keep the beat while the right tapper was instructed to synchronize/syncopate with the leader via a visual representation of the left tapper's tap on the screen.
    Right lead: identical to left lead except the rolls were reversed.
    Bidirectional: each tapper was given visual feedback from the other tapper and the two tappers were told to synchronize/syncopate with each other's tapping.}
    \label{fig:paradigms}
\end{figure}

For our analysis, we used data collected by Zhou, outlined in \cite{pinto2024symbolicdynamicsjointbrain}. 
We have included a short summary of the data collection process and a description of the data for the reader, but for more information, see the paper cited above.
The experiment consisted of six pairs of tappers, each pair consisting of a left and right tapper, isolated from each other. 
Each pair participated in a synchronized and syncopated session.
In the synchronized session, tappers tapped with the beat, while during the syncopated session, tappers tapped between the beat.
Each session consisted of 12 trials; four trial types (no lead, left lead, right lead, and bidirectional) repeated three times each where the trials were semi-randomly shuffled for each pair. 
Each tapper initially paced off of a metronome set at 1.3 Hz (78 beats per minute) for the first 30 beats and then tried to keep the beat after the metronome was removed. 
For additional explanation and visualization, see Figure~\ref{fig:paradigms}.

Each trial lasted $\sim3$ minutes and consisted of 230 taps. 
The first 30 taps of each trial, which used the metronome, were removed from the dataset.
Each session's dataset includes 32 channels of EEG data sampled at 2000 Hz and one channel that recorded the button presses for each participant. 
Additionally, the dataset included the session type, trial type, and channel labels.
We focused on differentiating between the synchronized and syncopated sessions for each tapper.

Just as in \cite{pinto2024symbolicdynamicsjointbrain}, we preprocessed the data prior to applying MapperEEG. 
Double taps were removed if two or more taps happen within 350ms. 
The taps were extracted from the pairs of subjects based on the subject with the smallest number of taps. 
Preprocessing the EEG data consisted of filtering out noise, detrending the data, re-referencing, applying a 0.5 Hz high-pass filter, and performing ICA to remove eye blinks.
Each session was split into trials and labeled by trial and session.
We downsampled each trial to reduce dataset size, computation time, and high frequency noise. 
Each synchronized and syncopated dataset contained over 3 million data points collected at a sampling rate of 2000 Hz.
We applied a cutoff frequency filter of 50 Hz and then downsampled each trial from 2000 Hz to 100 Hz, in accordance with Nyquist's theorem. 
This resulted in a 35-dimensional dataset consisting of 32 EEG channels, 1 response channel, 1 trial type channel, and 1 session type channel.

\subsection{Go/No-Go Shooting Task}\label{subsec:gonogo}

The second dataset consists of volunteers completing a go/no-go VR shooting task, introduced in Kerick et al.~\cite{kerick2023neural} and Spangler et al.~\cite{spangler2021gonogo}.
Each participant ideally completed six sessions, with some participants missing some sessions. 
A single session contained two difficulty conditions (low and high), and each condition consisted of four blocks with breaks between each block.
Each block contained 90 targets (trials) that the participants were instructed to shoot/not-shoot based on whether the target was an enemy (90\% of targets) or a friendly (10\% of targets).
These targets appeared randomly in one of nine established positions in the gaming environment.
During each session, participants completed 720 trials.

\begin{figure}
    \centering
    \includegraphics[width=0.75\linewidth]{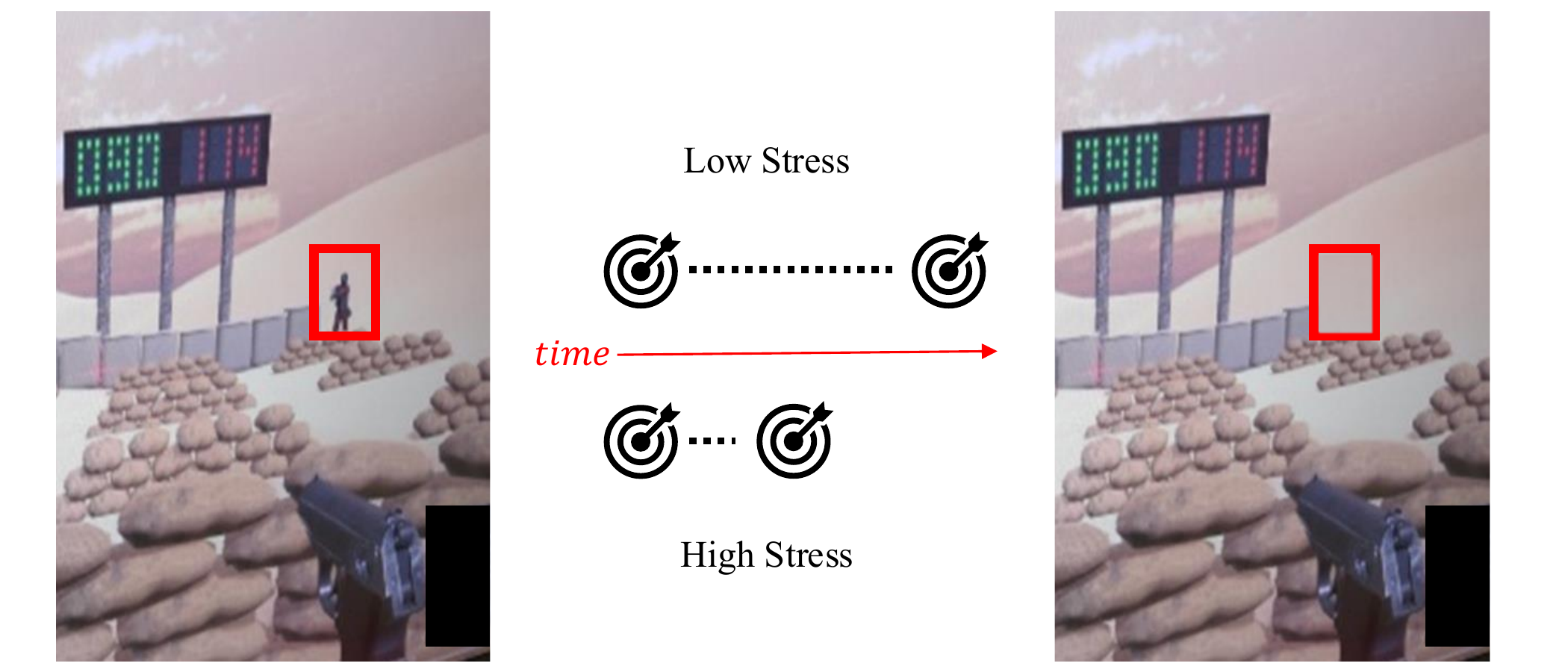}
    \caption{Screengrabs from the Go/No-Go Shooting task, with the left grab showing a target and the right showing no target. The low stress condition had longer target exposure time whereas the high stress condition had shorter target exposure time.}
    \label{fig:paradigms_go}
\end{figure}

Difficulty was thresholded based off of each individual. 
Participants had a chance to familiarize themselves with the equipment and task and were given time to practice before each session. 
As discussed in Spangler et al., ``Low difficulty was defined as the target exposure time (TET) that produced 90$\%$ shot accuracy on enemy targets, and high difficulty was defined as the TET that produced 50$\%$ shot accuracy on enemy targets'' \cite{spangler2021gonogo}.
After the low and high difficulty TETs were identified for each individual, trial TETs were randomly selected from corresponding Gaussian distributions centered on each participant's two average TETs, see Figure~\ref{fig:paradigms_go}.
Performance was calculated by a marksmanship value; $M = H/T$ where $H$ is the number of enemy targets hit and $T$ is the number of targets per block.

For our analysis we used data collected from 12 participants during their 6th session.
Session 6 is the last session the students participated in, and as such, they were familiar with the environment, instructions, and the tasks.
Researchers collected 64 channels of EEG data at 2048 Hz from each participant. 
Preprocessing included applying a ``high pass zero-phase Hamming-windowed sinc finite impulse response (FIR) filter with a cut-off frequency of 1 Hz'' \cite{Pratiher2022gonogo}.
Additionally, each channel was downsampled from 2048 Hz to 512 Hz. 
Researchers then applied the PREP preprocessing pipeline \cite{bigdely2015prep} and independent component analysis (ICA) \cite{jung2000removing} to remove artifacts. 
For additional information, see \cite{kerick2023neural, spangler2021gonogo, Pratiher2022gonogo}.
Preprocessing yielded a 65-dimensional dataset consisting of 64 EEG channels and 1 condition channel.

\section{Method: MapperEEG}\label{sec:Methods}

As a way to capture and cluster different brain states from EEG data, we look to the processes outlined in work such as those by Saggar, Geniesse, and others \cite{saggar2018towards, saggar2022precision, geniesse2019dyneusr}.
In each of these works, the researchers use Mapper as a way to characterize brain states based on fMRI data, but to our knowledge, this analysis has been restricted to fMRI data.
As such, we aim to show that similar structure can be captured and identified from EEG data.
In order to use Mapper as a way to cluster brain states captured by EEG data, we incorporate established EEG preprocessing methods with the Mapper algorithm and construct a clustering mechanism to categorize brain states and visualize brain state connectivity in a lower dimensional space.
Below, we outline the dataset and the framework of the MapperEEG algorithm.
In order to better establish intuition, we provide examples of how the steps are used on the tapping dataset and have framed the figure below in the same context.

\begin{figure}[h!]
    \centering
    \includegraphics[width=0.99\textwidth]{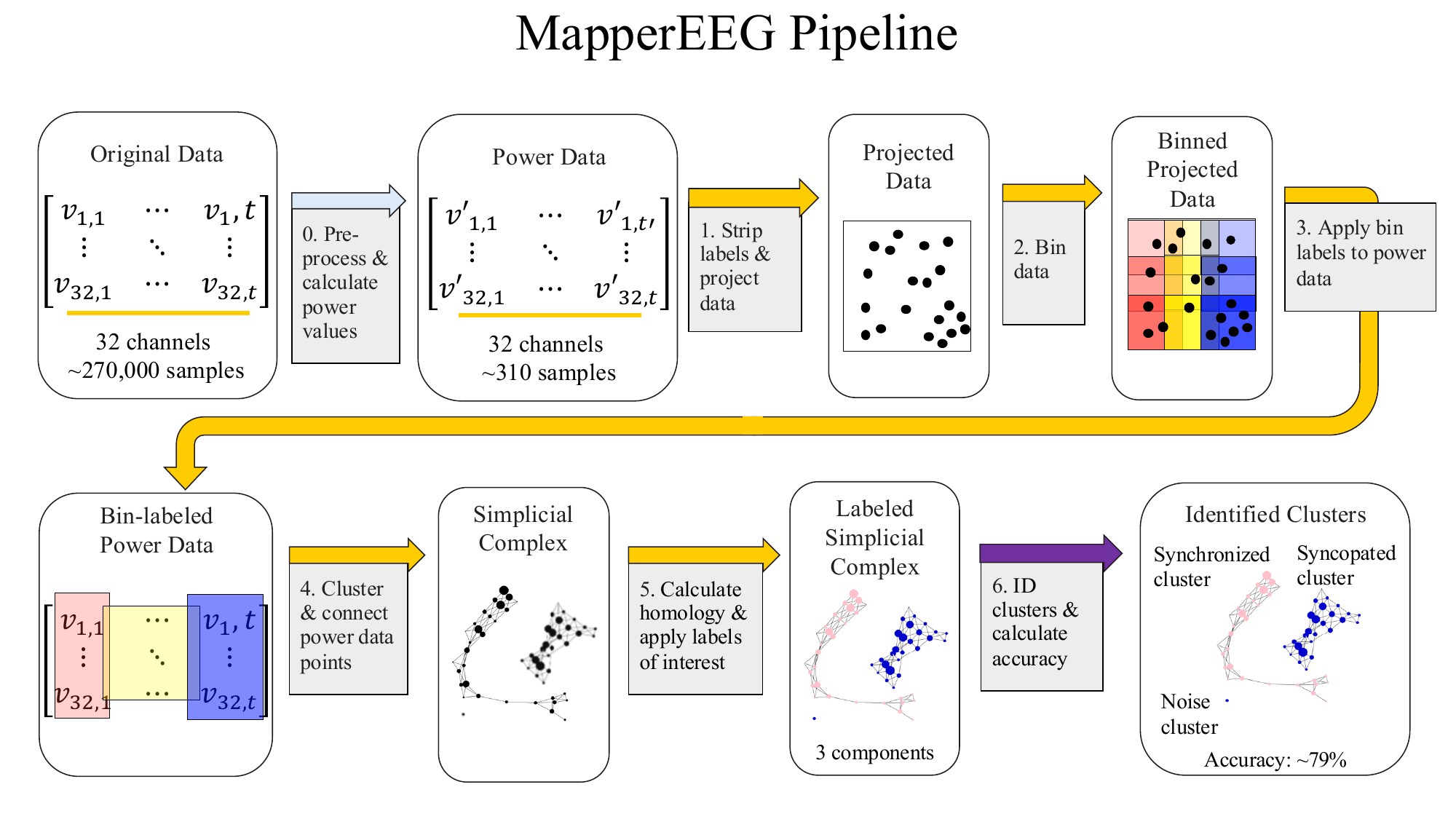}
    \caption{The MapperEEG pipeline. First (light blue arrow), we preprocess the data using well-established methods from EEG research. Next, we apply the Mapper algorithm (yellow arrows) as described in Subsection~\ref{subsec:Mapper}. Finally, we apply our clustering mechanism (dark purple arrow) that allows the Mapper algorithm to be used as a clustering tool. For additional information, see Section~\ref{sec:Methods}.}
    \label{fig:mapper}
\end{figure}

\subsection*{Step 0: Preprocess the data and construct the power dataset}\label{subsec:step0preprocess}

In our approach, we rely on established EEG techniques that capture frequency domain characteristics. 
To get a single power value per window, we calculated the amount of power in the traditional EEG bands: 1-4 Hz ($\delta$) 4-8 Hz ($\theta$), 8-13 Hz ($\alpha$) 13-30 Hz ($\beta$) 30-50 Hz ($\gamma$) as in \cite{Nayak2025eeg, vezard2015alertness}.
After the initial data cleaning and downsampling outlined in Section~\ref{sec:Datasets}, we calculated a sequence of power values by calculating and integrating the Power Spectral Density (PSD) functions for each channel using the welch python function from SciPy with a 1,000-sample window size and 50\% overlap.

\begin{figure}[ht]
    \centering
    \includegraphics[width=0.98\textwidth]{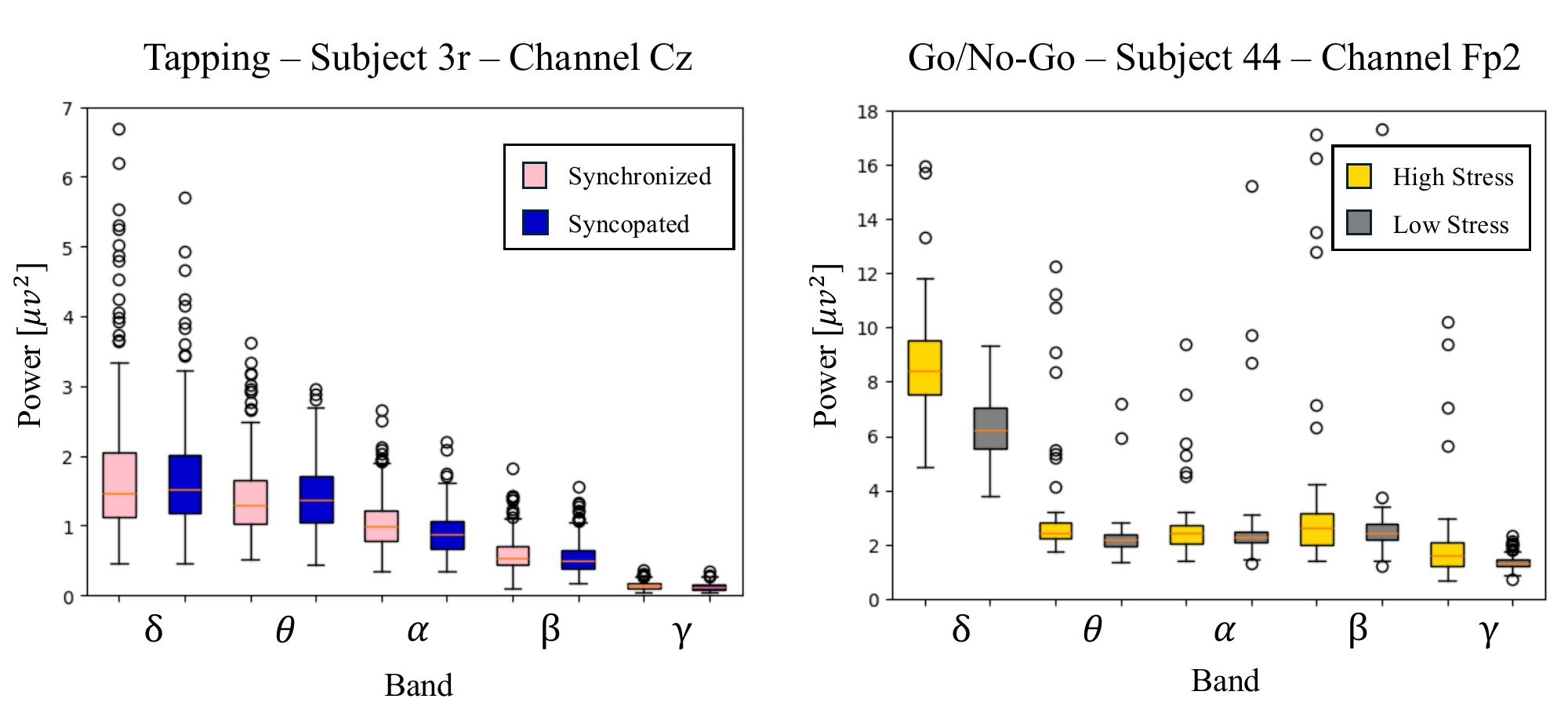}
    \caption{For illustration purposes, we show the distribution of power values for each band calculated for one subject in one channel for each of the two datasets.}
    \label{fig:power_box_plot}
\end{figure}

We restricted the upper bound of the gamma band from 80 Hz to 50 Hz and the lower bound of the alpha band from 0.1 Hz to 1 Hz to align with our downsampling and preprocessing algorithms. 
Note, other methods such as those found in Wang et al.~\cite{Wang2014-oo} and Belhadi et al.~\cite{Belhadi2025-gu}, show that using the frequency domain as opposed to the original time domain can lead to higher classification accuracy.
This gave us a sequence of power values for each of the five bands across each of the EEG channels. 
For an example of the resulting distribution of power values, see Figure~\ref{fig:power_box_plot}

\begin{figure}[ht]
    \centering
    \includegraphics[width=0.99\linewidth]{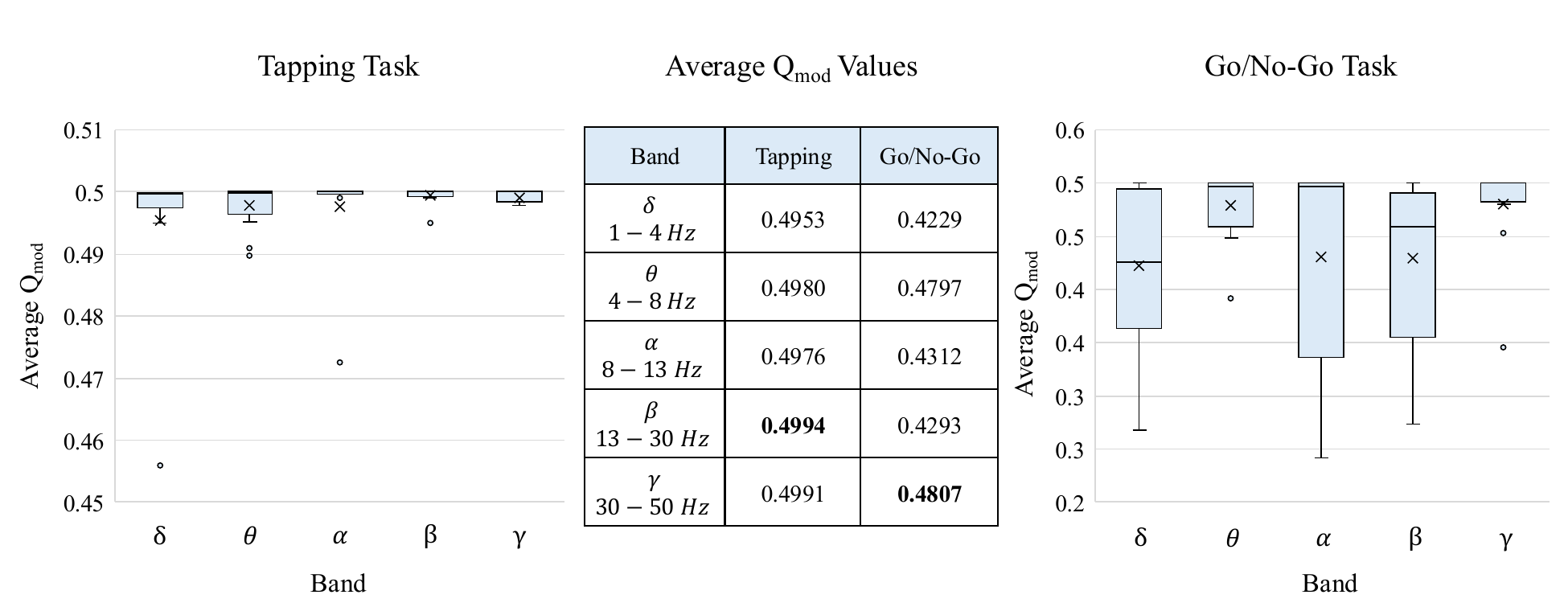}
    \caption{(Left and Right) Box and whisker plots capturing the average of the best \qmod values from each participant in each band. A single point on the box an whisker plot corresponds to the highest \qmod value for a single participant captured by running a grid search over different numbers of cubes and overlap amounts. (Center) The average \qmod value for each band and each task (corresponds to the $\times$'s on the box and whisker plots). Note, for the tapping task, all average values are within 0.005 of each other, thus we choose the gamma band based on the larger difference shown in the go/no-go task.}
    \label{fig:qmod_box}
\end{figure}

To identify which band contains the most salient features for unsupervised clustering, we applied MapperEEG to each of the five bands to identify the band that produced the highest average Q-mod score across both tasks via a hyperparameter search.
For both datasets, we applied MapperEEG to each participant's sessions across a range of parameters and calculated the \qmod score for each Mapper graph.
Recall, the \qmod score requires a label choice; for the tapping dataset, we used the synchronized/syncopated labels and for the shooting task, we used the high/low stress labels.
We performed a grid search across different numbers of bins (cubes in $\mathbb{R}^2$), $[10^2, 15^2, 20^2, 25^2, 30^2, 35^2, 40^2]$, overlap proportions, $[0.05, 0.1, 0.15, 0.2, 0.25, 0.3, 0.35, 0.4, 0.45, 0.5]$, and the five frequency bands, [1-4 Hz ($\delta$), 4-8 Hz ($\theta$), 8-13 Hz ($\alpha$), 13-30 Hz ($\beta$), 30-50 Hz ($\gamma$)].
For a summary of the results of the hyperparameter search for each band, see Figure~\ref{fig:qmod_box}.

All five bands produced an average \qmod score greater than 0.48 across all tappers and 0.25 across all shooting task participants, showing that MapperEEG is able to capture relevant information from all bands, see Table~\ref{fig:qmod_box}.
Ultimately, the gamma-band power produced the highest average \qmod score for each participant across all cube and overlap pairings in the go/no-go task, and as such, we chose the gamma-band power sequence as the Mapper algorithm input. 
For the tapping task, all average \qmod values are within 0.005 of each other, so we defer to the findings from the go/no-go task for choosing the optimal power band.

After preprocessing, downsampling, and converting the data into the spectral domain, we applied the Mapper algorithm to the gamma power value sequence, as first introduced in \cite{singh2007mapper}. 
The Mapper package used, Kepler Mapper, is outlined in \cite{hendrik2019keplermapper} and the code can be found in \cite{hendrik2020mappercode}.
As outlined above in Subsection~\ref{subsec:Mapper}, we followed a similar pipeline as that shown in \cite{saggar2018towards}.
In the following subsections, we describe each step of Mapper pipeline within the tapping experiment setting and how we modify the Mapper graph outputs to be used for clustering brain states.

\subsection*{Step 1: Project the data into a lower dimensional space}\label{subsec:step1project}

After generating the gamma power band dataset, we projected the power data from $\mathbb{R}^{32}$ dimensional space to $\mathbb{R}^2$ dimensional space such that $p:\mathbb{R}^{32}\rightarrow\mathbb{R}^2$. 
We tested three different projection functions: Principle Component Analysis (PCA) \cite{Jolliffe2002pca}, t-Stochastic Neighborhood Embedding (t-SNE) \cite{vander2008tsne}, and Uniform Manifold Approximation and Projection (UMAP) \cite{mcinnes2020umap}. 
PCA is a long-standing linear dimensionality reduction technique that transforms the data onto a new coordinate systems that captures the highest direction(s) of data variability.
The t-SNE algorithm is a statistical approach to data projection where data points that live close/far in the original space will be projected to points that live close/far in the low dimensional space with high probability.
UMAP is a non-linear dimension reduction algorithm based off manifold theory and topological data analysis \cite{mcinnes2020umap}. 
We chose t-SNE over UMAP and PCA due to it producing higher \qmod values in our preliminary testing compared to the other projection algorithms.

\subsection*{Step 2: Bin the projected data in the lower dimensional space}\label{subsec:step2bin}

Once the power data has been projected into $\mathbb{R}^2$, we cover the projected data with a square, $S\in\mathbb{R}^2$ such that $p(v)\in S$ for all $v\in \mathbb{R}^{32}$. 
Then, we tile the square with $b\times b$ square bins with $ov\%$ overlap. 
Note, choosing $b$ and $ov$ directly determines the side length of each bin.
Additionally, points can fall into multiple bins, but based on the construction of the square, each point will fall in at least one bin. 

\subsection*{Step 3: Apply lower dimensional bin labels to higher dimensional data}\label{subsec:step3label}

After each point has been assigned their bin label(s), we apply each point's bin label(s) to their pre-image points in the higher dimensional space. 
Note, one point might have multiple labels.

\subsection*{Step 4: Cluster the data into nodes and connect nodes}\label{subsec:step4cluster}

Next, we build the Mapper graph by clustering the points in the original, high dimensional space ($\mathbb{R}^{32}$) using the Density-Based Spatial Clustering of Applications with Noise (DBSCAN). 
Two points are clustered together into a node if:
\begin{enumerate}
    \item the two points share a bin label in $\mathbb{R}^2$, and
    \item if the two points are within a user-specified distance from each other in the high dimensional space ($\mathbb{R}^{32}$).
\end{enumerate}
Each node can contain multiple points and different nodes can contain the same point. 
In fact, two nodes are connected if and only if they share at least one data point, which is possible since points can be assigned multiple bin labels. 
The output is the Mapper graph, a graphical representation of how different brain states are grouped together in nodes, connected via edges, and how they evolve based on their lower-dimensional structure.

\subsection*{Step 5: Calculate the homology and apply labels of interest}\label{subsec:step5hom}

After building the simplicial complex, we applied simplicial homology and calculated the number of components of the Mapper graphs. 
Note, we could also calculate the number of holes and other higher dimensional features, but as we focused on clustering, we restricted our focus to the number of components. 

Next, we superimposed a pie chart denoting the label distribution of the data points contained in each node. 
The size of each pie chart was scaled relative to the number of points clustered into it's node, and the wedges of the chart corresponded to the distribution of labels within the node. 
Note, different labels can be used to capture different big pictures. 
As we will discuss in Section~\ref{sec:Discussion}, we evaluated multiple labeling schemes to get a broader understanding of the underlying structure of the Mapper graphs across two different datasets.

\subsection*{Step 6: Identify clusters and calculate accuracy}\label{subsec:step6ID}

To convert Mapper into a clustering algorithm, we identified each connected component as a cluster and assigned each of them a label. 
We assigned the label 1 to the cluster that contained the most points with label 1 (synchronized) and label 2 to the cluster with the most points with label 2 (syncopated) to ensure labeling consistency.
All additional clusters were grouped into a single noise cluster, meaning the algorithm did not cluster them correctly, such as how some points are assigned to the noise cluster in the DBSCAN algorithm.

Note, each data point can only be assigned to a single cluster.
If one data point falls into two different nodes, those nodes are connected with an edge, which merges the two components.
Thus, the only way an element can be in multiple graph nodes is if those nodes are each in the same component.

Finally, we calculated the cluster accuracy using five different metrics: F-1 Average, F-1 Macro, F-1 Weighted, Silhouette, and Davies-Bouldin.

\section{Results}\label{sec:Results}

Below, we compare MapperEEG to six other unsupervised clustering algorithms on the gamma-band power sequences generated from each of the EEG datasets. 
Specifically, we compare MapperEEG against DBSCAN, $k$-Means, Gaussian Mixture Model, Hierarchical Clustering, an Autoencoder, and a Hidden Markov Model using the F-1, Silhouette, and Davies-Bouldin metrics.
All experiments were performed on a 2023 MacBook Pro with an M3 Max Chip and 64 GB of memory.
The average computational complexity of MapperEEG is dependent on the computational complexity of the projection function and the clustering function. 
For UMAP, the average computational complexity is $O(N^{1.14})$ 
and DBSCAN is $O(N\cdot \log N)$ \cite{mcinnes2020umap, Ester1996DBSCAN}.
Thus, we estimate the average computational complexity of MapperEEG to be $O(N\cdot \log N + N^{1.14})$.
Below, we compare MapperEEG against the other clustering algorithms on the tapping and go/no-go datasets.

\subsection{Tapping Dataset}\label{subsec:compare_tap}

We analyzed how MapperEEG performs as an unsupervised clustering algorithm by evaluating how well it identifies states correctly as synchronized or syncopated for the tapping task.
In Subsection~\ref{subsec:step0preprocess}, we identified the gamma band as the best band choice for the MapperEEG algorithm based on its resulting high community structure. 
We tested MapperEEG against six common clustering algorithms across five different error measures, F-1 Micro (accuracy), F-1 Macro, F-1 Weighted, Silhouette Score, and Davies-Bouldin. 
For each algorithm, we performed a hyperparameter search across a range of values, and identified the hyperparameters that yielded the highest average accuracy across all participants. 
For DBSCAN, we allowed the minimum number of samples ($ms$) to range by 4's between $[4, 64]$ and the maximum distance between samples of the same cluster ($eps$) range between $[0.5, 10]$ by steps of $0.5$.
For $k$-Means and for hierarchical clustering we allowed the number of clusters to vary across the interval $[2,10]$.
Similarly, for the Gaussian Mixture Model (GMM) and the Hidden Markov Model (HMM), we looked across $[2,10]$ clusters/components and across four coverage types, [full, tied, diag, sphere]. 
To test the autoencoder, we varied the number and dimension of the layers in both the encoder and the decoder, $\big[[32-16-8-4-2-4-8-16-32], [32-16-4-2-4-16-32], [32-16-2-16-32]\big]$ and the number of training epochs over $[20, 40, 60, 80, 100]$. 

The results for the parameters that yield the best average accuracy error across all participants for the tapping experiment are included in Table~\ref{tab:errors_tapping}.

\begin{table}[h]
    \centering
    \begin{tabular}{| c | c || c | c | c | c | c |}
    \hline 
        Alg/Metric & Parameters & Accuracy & F-1 Macro & F-1 Wgt & Silhouette & DB\\ 
        \hline\hline
        DBSCAN  & $ms=12$, $ep=7.5$   & 0.47889 & 0.25442 & 0.32457 & 0.56480 & \textbf{0.54817}\\
        \hline
        k-Means & $clus=2$            & 0.58217 & 0.35087 & 0.52556 & 0.56457 & 1.19653\\
        \hline
        Hierarchical & $clus=4$       & 0.63425 & 0.37167 & 0.55669 & \textbf{0.59673} & 0.95960\\
        \hline
        GMM     & $cps=5$, $var=full$ & 0.72928 & 0.47949 & 0.71864 & 0.43484 & 1.43536\\
        \hline
        HMM     & $cps=4$, $var=full$ & 0.74348 & 0.48411 & 0.72562 & 0.40548 & 1.64403\\
        \hline
        AE      & layers=8, epochs=40 & 0.63937 & \textbf{0.52025} & 0.57267 & 0.65463 & 0.62434 \\
        \hline
        MapperEEG & cbs=20, ovlp=35\% & \textbf{0.82497}& 0.30074 & \textbf{0.85666} & -0.02709 &  1.52661\\
        \hline
    \end{tabular}
    \caption{The average metrics for each of the seven clustering algorithms. Best scores are bolded for ease of comparison. For all metrics other than DB, higher is better. For DB, lower is better.}
    \label{tab:errors_tapping}
\end{table}

MapperEEG outperformed the other algorithms across the Accuracy and F-1 Weighted error functions. 
These errors are calculated by comparing the true label list to the predicted label list and identifying how accurate each unsupervised clustering algorithm is.
The last two scores capture how well the algorithm preserves the distance between points.
The silhouette score calculates the distance between a sample and the nearest cluster that the sample is not a part of.
The Davies-Bouldin score evaluates the ratio of within-cluster distances to between-cluster distances which means clusters that are internally less dispersed and farther apart from each other will result in a better score.
Both of these metrics evaluate how the algorithm clusters data in the original domain ($\mathbb{R}^{32}$).
Since MapperEEG projects the data and calculates the labels based off of those projections, MapperEEG is capturing structure that is not strictly dependent on how close the points are in the original domain space. 
So even though MapperEEG has a worse Silhouette and Davies-Bouldin score than some of the other methods, MapperEEG's clustering accuracy is significantly higher. 
This shows that the MapperEEG algorithm is capturing structure that is not easily identified in the higher dimensional space.  

\subsection{Go/No-Go Dataset}\label{subsec:compare_gonogo}

In order to show that MapperEEG usefulness is not restricted to a single dataset, we performed the same evaluation on the gamma power data from the go/no-go dataset. 
We compared MapperEEG to the same six clustering algorithms using the same five error measurements; DBSCAN, $k$-Means, Gaussian Mixture Method (GMM), Hierarchical Clustering, an autoencoder (AE), and a Hidden Markov Model (HMM) using the three F-1 metrics, the Silhouette score, and the Davies-Bouldin (DB) metric.
Additionally, we used the same hyperparameter options as outlined above in Subsection~\ref{subsec:compare_tap}.
The results for the parameters that yield the best average accuracy error across all participants for the go/no-go experiment are included in Table~\ref{tab:errors_gonogo}.

\begin{table}[h]
    \centering
    \begin{tabular}{| c | c || c | c | c | c | c |}
    \hline 
        Alg/Metric & Parameters & Accuracy & F-1 Macro & F-1 Wgt & Silhouette & DB\\ 
        \hline\hline
        DBSCAN  & $ms=7$, $ep=5$                & 0.41879 & 0.21187 & 0.27919 & 0.75447 & 0.95748\\
        \hline
        k-Means & $clus=2$                      & 0.47703 & 0.26194 & 0.33134 & \textbf{0.87797} & 0.43540\\
        \hline
        Hierarchical & $clus=2$                 & \textbf{0.52784} & \textbf{0.40566} & 0.39989 & 0.87464 & 0.42825 \\
        \hline
        GMM     & $cps=10$, $var=sph$   & 0.48247 & 0.33192 & 0.48899 & 0.34276 & 1.55120\\
        \hline
        HMM     & $cps=10$, $var=sph$   & 0.50245 & 0.34332 & \textbf{0.50448} & 0.37259 & 1.44243\\
        \hline
        AE      & layers=6, epochs=60           & 0.52560 & 0.35043 & 0.40872 & 0.83517 & \textbf{0.42668}\\
        \hline
        MapperEEG & cbs=10, ovlp=25\%           & 0.41250 & 0.11584 & 0.32870 & 0.03517 & 1.57885 \\
        \hline
    \end{tabular}
    \caption{The average metrics for each of the seven clustering algorithms. Best scores are bolded for ease of comparison. For all metrics other than DB, higher is better. For DB, lower is better.}
    \label{tab:errors_gonogo}
\end{table}

For this dataset, all methods fail to differentiate between the high stress and low stress conditions. 
Further, MapperEEG does not outperform the other algorithms, and instead, falls below some of the other methods. 
But, in Subsection~\ref{subsec:step0preprocess}, we showed that the MapperEEG outputs have high community structure across both the tapping and the go/no-go shooting tasks. 
Although the community score for the tapping task is higher, the go/no-go task still yields significant community structure. 
Thus, we turn to the Mapper graphs to better understand the underlying structure of the data.

\subsection{Beyond Clustering}\label{subsec:beyondclustering}

The ability to capture and understand brain state connectivity is critical for providing researchers with insight into how brain states evolve over time.
Mapper graphs provide a representation of that evolutionary structure. 
In Figure~\ref{fig:spec_v_temporal}, we show how there is high connectivity between similar states, even if they aren't partitioned into different clusters.
The high community structure across most participants demonstrates that even without distinct clusters, MapperEEG can provide additional insight into the underlying groupings of brain states.

\section{Discussion}\label{sec:Discussion}

This work focused on introducing the MapperEEG method as a topological approach to analyzing EEG data, specifically it's ability to cluster brain states in a structurally relevant way. 
We focused this introduction on identifying the differences between individual brain states characterized by EEG recordings between two experiments with two different experimental conditions; synchronized and syncopated finger tapping and high stress and low stress go/no-go tasks.
Below, we outline the importance of the spectral domain, current limitations, and MapperEEG beyond clustering.

\subsection{The Importance of the Spectral Domain}\label{subsec:spectral}

Since EEG's foundation, researchers have used different techniques to transform data from the time domain to the spectral domain as a way to understand and find patterns in the data \cite{ng2022spec}.
Unlike other data formats, such as fMRI, EEG data in the temporal domain is dominated by noise which makes it harder to identify patterns and other salient information. 
Thus, we show the efficacy of switching to the spectral domain on the tapping data and the go/no-go data described above. 

\begin{figure}[ht]
    \centering
    \includegraphics[width=\textwidth]{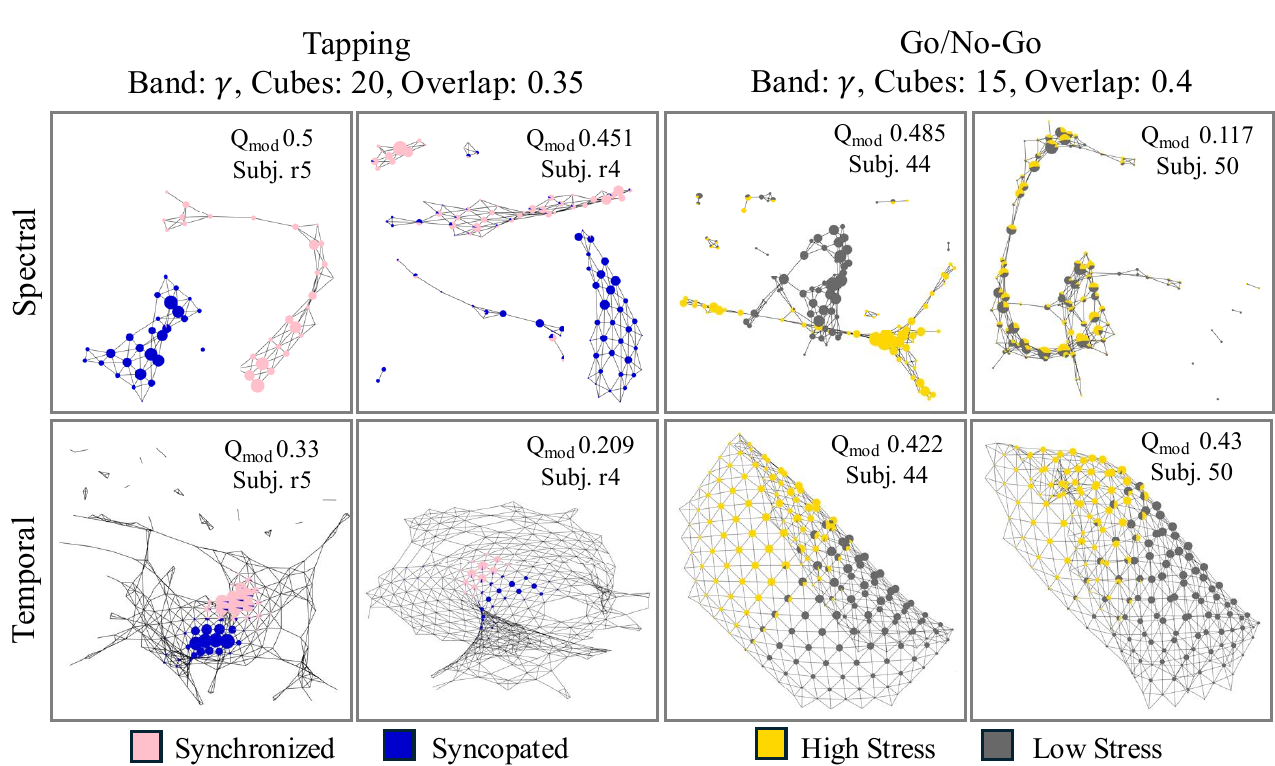}
    \caption{Examples of Mapper graphs for the tapping (columns 1 and 2) and go/no-go (columns 3 and 4) datasets in the spectral and temporal domains. Each column's Mapper graph is generated from one subject's EEG data. In the Spectral row, the data was preprocessed into the gamma power dataset whereas in the Temporal row, the data was left as is. For each experiment, the MapperEEG parameters were fixed and subjects for illustration were chosen based on those that had the highest spectral \qmod scores (columns 1 and 3) and the lowest spectral \qmod scores (columns 2 and 4).}
    \label{fig:spec_v_temporal}
\end{figure}

As seen in first two columns of Figure~\ref{fig:spec_v_temporal}, applying Mapper to the temporal tapping data yields a overly-connected Mapper graph with less organization than the spectral Mapper graphs.
Further, the temporal tapping Mapper graphs have lower \qmod scores, showing there is less community structure in these graphs.
Additionally, for the spectral graphs in the tapping dataset, even the Mapper graph with the lowest \qmod score still contains two distinct clusters that partition the synchronized and syncopated states.
With the overlayed synchronization and syncopation labels, it is clear that Mapper is able to cluster the spectral data appropriately, whereas it was unable to do as well with the temporal data.

Comparatively, the Mapper graphs for the go/no-go task show different findings. For subject 44, converting the data to the spectral domain led to distinct clusters within the Mapper graphs (column 3, row 1). 
But for subject 50, there were not any clear clusters within the Mapper graph (column 4, row 1). 
The \qmod score reflects this lack of community.
But, if we look at the temporal Mapper graphs for the go/no-go tasks, even though there are not any clusters, the data is still partitioned into 2 parts of the Mapper graph.
The \qmod scores support this observation, and these scores are close to those of the spectral tapping Mapper graphs, even if the brain states are not partitioned into distinct clusters. 
Thus, even though our results showed that as clustering algorithm, MapperEEG did not outperform it's competitors, it was still able to capture the underlying structure of the data in the temporal domain via graph connectivity.

This finding is not terribly surprising. 
The existing body of work for EEG analysis shows that different brain states are correlated with the presence of particular brain rhythms, for example, see Buzsaki et al.~\cite{buzsaki2006rhythms}. 
Thus, transforming the EEG signal into the spectral domain can provide direct information about the brain state rhythms that might be lost in the noise of the temporal domain.
Specifically for finger tapping experiments, researchers have shown that spectral analysis of EEG data can be used to identify finger movements~\cite{Xiao2015finger}. 
Hence, it is not too surprising that switching to the frequency domain provided better input data for the MapperEEG algorithm than the time domain data on the tapping dataset. 
In another experiment where time domain data would be more relevant, such as the go/no-go task, it is feasible MapperEEG would work better on the temporal EEG data. 
However, we expect that since the existing body of work in EEG analysis provides evidence that particular rhythms in the brain are relevant, it is prudent to investigate if the frequency domain would provide clearer results than the time domain. 
Finally, there is also possibility that applying another preprocessing approach before applying any classification algorithms could be useful. For example, calculating complexity measures of the time series \cite{kerick2023neural}.

\subsection{Tapping versus Go/No-Go: Why gamma works for one and not the other}\label{subsec:differences}

As shown in Figure~\ref{fig:spec_v_temporal}, gamma-band features work well for distinguishing between the two tapping tasks, while they prove less discriminatory for the go/no-go task. 
One reason this could be is that gamma-band oscillations are typically associated with local cortical processing and precise temporal coordination of neuronal populations~\cite{Fries2015gamma}.  
During motor tapping tasks, which require rhythmic execution and often involve neural entrainment, gamma activity reflects tight local synchrony within sensorimotor networks—essentially mediating fine-grained temporal coordination between perception and action that differs between synchronized and syncopated states. 
This may explain why gamma features are more discriminative in the motor synchronization condition.

In contrast, during cognitive stress or inhibitory control tasks (such as the go/no-go task), neural dynamics typically engage broader, distributed networks (e.g., frontal–parietal control circuits) and greater reliance on theta and beta rhythms related to cognitive monitoring and top-down control~\cite{Engel2010stress}. 
Gamma activity may still occur locally, but it tends to be less spatially coherent and less consistent across participants, which could account for the weaker separation observed in the Mapper graphs for that task.
Note, this is a theoretical interpretation based on established neurophysiological frameworks; functional connectivity between brain regions was not directly examined in either study.

\subsection{Limitations}\label{subsec:limitations}

As touched on throughout, MapperEEG does have its limitations.
MapperEEG's computational complexity ($O(N\cdot \log N + N^{1.14})$) is greater than several of the algorithms we compare against, such as DBSCAN ($O(N\cdot \log N)$ \cite{Ester1996DBSCAN}), $k$-Means ($O(NKDI)$), and the Gaussian Mixture Method ($O(NKD^3$).
MapperEEG's complexity can be reduced by changing the projection and clustering algorithms to ones with lower complexity, such as changing the projection function to PCA.

Additionally, MapperEEG is not effective in the temporal domain, unlike other clustering algorithms such as microstate analysis or HMM~\cite{Jajcay2023micro, si2021mirco}. 
Applying the Mapper algorithm directly to the EEG time series did not separate brain states during the different conditions of the experiment. 
Only after performing spectral analysis of the EEG data and using power within different EEG spectral bands could we successfully cluster data into groups associated with the syncopated and synchronized conditions. 

There is also relatively high sensitivity to the choice of the number of cubes and the percent overlap.
This can be accounted for by performing a hyperparameter sweep, but that takes time and the best hyperparameters are not always consistent across participants.
Further, to fully investigate all potential hyperparameter combinations, one would also have to sweep over different values for the PSD window size, the projection algorithm hyperparameters, and the clustering algorithm hyperparameters, as all affect the Mapper graph output. 

\section{Conclusions and Future Work}\label{sec:Conclusion}

Capturing the underlying structure of brain states is an ongoing endeavor. 
Tools that allow us to visualize connectivity and structure without relying on known labels provide researchers with an insight into brain state relationships.
Methods that allow for visualization of this structure will help identify individuals who have high compartmentalization levels (high modularity) and allow us to compare the underlying brain structure across individuals. 
For example, some participant's Mapper graphs had higher average \qmod scores across different MapperEEG parameters, such as Tappers 1, 5, and 9 whereas others, such as Tapper 10 demonstrate significantly lower average \qmod scores (see Figure~\ref{fig:qmod_box}). 
Finding that some of these individuals have more community structure than others might allow us to better understand why some people do better at some tasks than others. 
In the work above, we demonstrate that MapperEEG produces brain state connectivity graphs with significant community structure by producing a high \qmod value across a range of Mapper values.
Additionally, MapperEEG demonstrated its effectiveness as a clustering algorithm as it outperformed DBSCAN, $k$-Means, and the Gaussian Mixture Method by over 20\% with respect to accuracy.
These two pieces show that MapperEEG has the potential to be a useful structural assessment and visualization tool for analyzing EEG data.

For this analysis, we restricted our focus to the \qmod metric and how it could be used to choose MapperEEG parameters across a wide range of hyperparameters.
This is just the beginning of how the Mapper graphs can be used to develop a deeper understanding of brain states and their functional connectivity.
The results in \cite{saggar2018towards} show how Mapper is able to capture transition paths between brain states. 
In the same work, Saggar et al. also looked at how core-periphery structure corresponds to task performance. 
Our future work will look at investigating whether the same kinds of transitions and underlying structure found within fMRI data can be captured with EEG.
Specifically, we want to explore how nodes are connected with respect to time to see how states transition across time.

The datasets we used for this initial paper did not record brain activity between trials, which limited the relevant analyses that could be used.
Thus, for future projects, we will to design experiments that record EEG data during and between tasks.
This will allow us to use more tools from network science, topology, and graph theory to learn more about the underlying functional connectivity of the brain.
Additionally, we plan to separately analyze data from individual EEG channels and look at how source localization might provide more information on brain state.

There are also other preprocessing methods we could look into, such as computing the complexity values of the signal~\cite{kerick2023neural} or using information about the aperiodic and periodic activity within the PSD~\cite{yu2024aperiodic}.
Moving forward, we plan to continue to use MapperEEG to explore brain state connectivity and we hope to use that knowledge to understand and/or predict an individual's brain state during a task.
This would require MapperEEG to be used in real-time to monitor patient brain state transitions and provide their current brain state. 
Real-time monitoring requires an algorithm that is fast, adaptive, and redesigned to reduce the computation resource requirements. 
The current implementation of MapperEEG takes approximately $\sim 16.8$ seconds to run the current MapperEEG pipeline on a dataset of $\sim$270,000 data points.  
As the number of data streams, amount of data, and analytical approaches increases, we plan to investigate GPU implementations of Mapper, such as those in Zhou et al.~\cite{zhou2020mapper}, to increase the efficiency of MapperEEG.
Further experiments will be designed to investigate how the Mapper graph changes when new data is introduced.
Will the graph maintain it's underlying structure, or will it evolve over time?
How much past data is needed in order to preserve the underlying structure while allowing for new data to be incorporated in?
These questions are critical to exploring how MapperEEG can go beyond the applications explored in this paper.

\bibliographystyle{unsrt}
\bibliography{main}

\end{document}